\documentclass[12pt, twoside]{article}
\usepackage{amsmath, amssymb, amsthm}
\usepackage{enumerate}
\usepackage{tikz} 
\usepackage{natbib}

\usepackage{graphicx,subcaption}
\usepackage{array}

\usepackage{graphicx,color}

\newcommand{\coloredcircle}[1]{\tikz{\draw[fill=#1,draw=#1] (0,0) circle (0.15cm);}} 

\theoremstyle{plain}
\newtheorem{theorem}{Theorem}[section]

\newtheorem{lemma}{Lemma}[section]

\newtheoremstyle{boldremark}
{\topsep}
  {\topsep}
  {\normalfont} 
  {} 
  {\bfseries} 
  {.} 
  { } 
  {} 
\theoremstyle{boldremark} 
\newtheorem{remark}{Remark}[section] 

\def\oo{{\bf 0}}

\def\oo{{\rm {o}}} 
\def\OO{{\rm {O}}}

\begin{document}

\baselineskip 6true mm
	
	\bigskip
	\centerline{\Large\bf Limit theorems for runs }
    \centerline{\Large\bf containing two types of contaminations.}
	\centerline{\Large\bf Paper with detailed proofs}
	\vskip0.3true cm
	\centerline{
		\renewcommand{\thefootnote}{1}
		Istv\'an Fazekas\footnotemark , \ Borb\'ala Fazekas, Michael Ochieng Suja
		\footnotetext{fazekas.istvan@inf.unideb.hu}
}
\vskip0.3truecm
\centerline{University of Debrecen, Hungary}
\vskip0.4truecm

\begin{abstract}
In this paper, sequences of trials having three outcomes are studied.
The outcomes are labelled as success, failure of type I and failure of type II.
A run is called at most $1+1$ contaminated, if it contains at most $1$ failure of type I
and at most $1$ failure of type II.
The limiting distribution of the first hitting time 
and the accompanying distribution for the length of the longest at most $1+1$ contaminated run are obtained.
This paper contains the detailed mathematical proofs.
Simulation results supporting the theorems are also presented.
\end{abstract}
  \section{Introduction}
 The problem of the length of the longest head run for $n$ Bernoulli random variables was first raised by T. Varga. 
 The first answer for the case of a fair coin was given in the classical paper by \citet{erdos1970law}.
Surprisingly more precise answer, the almost sure limit results for the length of 
the longest runs containing at most $T$ tails was given by \citet{erdos1975length}.
Their result is the following.
Consider the usual coin tossing experiment with a fair coin.
Let $Z_T(N)$ denote the longest head run containing at most $T$ tails in the first $N$ trials.
Let $\log$ denote the logarithm to base $2$ and let $[ . ]$ denote the integer part.
Let $h(N)= \log N +T \log \log N- \log \log \log N -\log T! + \log \log e$.
Let $\varepsilon$ be an arbitrary positive number.
Then for almost all $\omega\in \Omega$ there exists a finite $N_0=N_0(\omega)$
such that $Z_T(N)\ge [h(N) -2 -\varepsilon]$ if $N\ge N_0$,
moreover, there exists an infinite sequence $N_i=N_i(\omega)$ of integers such that
$Z_T(N_i) < [h(N_i) -1 +\varepsilon]$.
 
These results later inspired renewed interest in this research area and several subsequent papers came up. 
Asymptotic results for the distribution of the number of $T$-contaminated head runs, the first hitting time of a 
$T$-contaminated head run having a fixed length, 
and the length of the longest $T$-contaminated head run were presented by \citet{foldes1979limit}.
For the asymptotic distribution of $Z_T(N)$, \citet{foldes1979limit} presented the following result.
For any integer $k$ and $T\ge 0$ 
$$
P(Z_T(N) - [\log N + T \log \log N ] < k ) = \exp\left(- \frac{2^{-(k+1-\{\log N +T \log \log N \} )} }{T!} \right)
+ \oo(1) 
$$
as $N\to \infty$, where $\{ . \}$ denotes the fractional part.

By applying extreme value theory, \citet{gordon1986extreme} obtained the asymptotic behaviour of 
the expectation and the variance of the length of the longest $T$-contaminated head run. 
Also in the same, \citet{gordon1986extreme}, accompanying distributions were obtained for the length of the
longest $T$-contaminated head run.

\citet{fazekas2021limit} shown that the accompanying distribution initially obtained by \citet{gordon1986extreme} 
could as well be arrived at using the approach given by \citet{foldes1979limit}.
After some probabilistic calculations and algebraic manipulations of the approximation of the length of the 
longest $T$-contaminated run and using the main lemma in \citet{csaki1987limit}, 
more precise and accurate convergence rate was obtained for the accompanying distribution of 
the longest $T$-contaminated head run by \citet{fazekas2023convergence} for $T=1$ and $T=2$.

Different authors have given in depth considerations to experiments involving sequences of runs emerging from 
trinary trials, where Markov chain approach is used in their analysis. 
Such sequences include system theoretic applications, where components might exist in the following states; 
`perfect functioning', `partial functioning' and `complete failure' (see \cite{eryilmaz2016generalized}
 and \cite{koutras1997sooner}).

In this paper, we study sequences of trials having three outcomes: success, failure of type I and failure of type II.
We shall say that a run is at most $1+1$ contaminated, if it includes at most $1$ failure of type I
and at most $1$ failure of type II.
We shall study the limiting distribution of the first hitting time 
and the accompanying distribution for the length of the longest at most $1+1$ contaminated run.
Section \ref{results} contains the main results.
It is proved that the appropriately normalized first hitting time has exponential limiting distribution, 
see Theorem \ref{T1}. 
Then we give accompanying distributions for the appropriately centralized length of the longest at most 
$1+1$ contaminated run, see Theorem \ref{T2}.
The proofs are presented in Section \ref{proofs}.
The proofs are based on a powerful lemma by \citet{csaki1987limit}.
For reader's convenience, we quote the lemma in the Appendix section \ref{app}.
In this paper we describe the details of the quite long proofs.
In Section \ref{simu} we present simulation results supporting our theorems.
 \section{The longest at most $1+1$ contaminated run and the first hitting time}  \label{results}
 Let $X_1,X_2,\ldots,X_N$ be a sequence of independent random variables with three possible outcomes; 
 $0$, $+1$ and $-1$ labelled as success, failure of type I and failure of type II, respectively with the distribution
 $$
 P(X_i = 0) = p,  \quad P(X_i= +1) = q_1 \ \  {\text {and}} \ \ P(X_i= -1) = q_2,
 $$ 
 where $p+ q_1+q_2 = 1$ and $p>0$,\quad$q_1>0$,\quad$q_2 > 0$.
 
 An $m$ length sequence is called a pure run if it contains only $0$ values. 
 It is called a one-type contaminated run if it contains precisely one non-zero element either a $+1$ or a $-1$. 
 On the other hand, it is called a two-type contaminated run if it contains precisely one $+1$, and one $-1$
  while the rest of the elements are $0$'s.
  
 A run is called at most one+one contaminated (shortly at most $1+1$ contaminated), 
 if it is either pure, or one-type contaminated, or two-type contaminated.

 So for an arbitrary fixed $m$, let $A_n = A_{n,m}$ denote the occurrence of the event at the $n^{th}$ step, 
 that is, there is an at most $1+1$ contaminated run in the sequence $X_n, X_{n+1}, \ldots, X_{n+m-1}$ 
 and $\bar{A}_n$ be its non-occurrence.
 We see that
 $$
 P(A_1) = p^m + m(1-p)p^{m-1} + m(m-1)p^{m-2} q_1 q_2.
 $$

In what follows, we shall use the notation
 \begin{equation*}
      \alpha = 
      \frac{C_0 +\frac{1}{m} C_1 +\frac{1}{m(m-1)}C_2  }{1+ \frac{p(1-p)}{(m-1)q_1q_2} + \frac{p^2}{m(m-1)q_1q_2}} ,
  \end{equation*}
  where
  $$
  C_0=(q_1+q_2), \ C_1= \frac{p(q_1^2+q_2^2)}{q_1q_2}-1, \ 
  C_2= \frac{(q_1^2+q_2^2)p^2}{q_1q_2(p-1)} + \frac{p}{p-1} + \frac{2(2p+1)q_1q_2}{(p-1)^3}.
 $$
 \subsection{The first hitting time of the at most $1+1$ contaminated run of given length}
 Let $\tau_m$ be the first hitting time of the at most $1+1$ contaminated run of heads having length $m$. 
 We shall be interested in finding the limiting distribution of $\tau_m$ as $m \to \infty$ 
 for the case of a sequence containing at most two contaminations but no two of the same type.
\begin{theorem}\label{T1}
 Let $P(X_i = 0) = p$, $P(X_i= +1) = q_1$ and $P(X_i= -1) = q_2$ be probabilities of success, 
 failure of type I and failure of type II, respectively, 
 where $p+ q_1+q_2 = 1$ and $p>0$, \ $q_1>0$, \ $q_2 > 0$. 
 Let $\tau_m$ be the first hitting time of the at most $1+1$ contaminated run of heads having length $m$. 
 Then, for $x>0$,
 \begin{equation}
     P(\tau_m \alpha P(A_1)> x) \sim e^{-x}
 \end{equation}
 as $m\to \infty$. 
\end{theorem}
 \subsection{Length of the longest at most $1+1$ contaminated run}
Let $\mu(N)$ be the length of the longest at most $1+1$ contaminated run in $X_1,X_2,\ldots,X_N$. 
Then,
 $\{\mu(N) < m\}$ if and only if any $m$ length run in $X_1,X_2,\ldots,X_N$ 
 is neither two-type contaminated nor one-type contaminated nor pure.
 
 We need some notation.
Let 
$$
K = \frac{2C_0C_2 - C_1^2 - C_0^2}{2CC_0^2},
$$
where $C=\ln\frac{1}{p}$ and $\ln $ is the logarithm to base $e$.
Let
\begin{equation*}
    \begin{split}
 {}&   m(N) =
    \log N + 2 \log \log N + \frac{4 \log \log N}{C \log N} + \frac{C_1-C_0}{C C_0} \frac{1}{\log N} - \\
    & -\frac{4}{C} \frac{(\log \log N)^2}{(\log N)^2} + \left( \frac{8}{C^2} - 
     \frac{2(C_1-C_0)}{CC_0} \right)\frac{\log \log N}{(\log N)^2} + 
     \left(\frac{2(C_1-C_0)}{C^2 C_0} + K\right)\frac{1}{(\log N)^2} + \\
     &+  \frac{16}{3C} \frac{(\log \log N)^3}{(\log N)^3} + \left( -\frac{16}{C^2} + \frac{4(C_1-C_0)}{C C_0}\right)
     \frac{(\log \log N)^2}{(\log N)^3} - 
     \left(4K + \frac{8(C_1-C_0)}{C^2 C_0} \right)\frac{\log \log N}{(\log N)^3} + \\
    & + \frac{16}{C^3}\frac{\log\log N}{(\log N)^3} - \frac{8}{C^2}\frac{(\log\log N)^2}{(\log N)^3} - 
        \frac{4(C_1-C_0)}{C^2 C_0} \frac{\log \log N}{(\log N)^3} ,
     \end{split}
\end{equation*}
where
$\log$ denotes the logarithm to base $1/p$.   
Let $[m(N)]$ denote the integer part of $m(N)$ and let $\{m(N)\}= m(N)-[m(N)]$ denote its fractional part.
 Introduce the function
\begin{equation*}
    \begin{split}
        H(x) {}& =
        -x + \frac{2}{C\log N}x - \frac{4}{C}\frac{\log \log N}{(\log N)^2}x -
        \frac{C_1-C_0}{C C_0} \frac{1}{(\log N)^2}x + \\
        + {}&\left( \frac{4(C_1-C_0)}{C C_0} - \frac{8}{C^2}\right)\frac{\log \log N}{(\log N)^3}x + 
        \frac{8}{C}\frac{(\log \log N)^2}{(\log N)^3}x - \frac{1}{C}\frac{1}{(\log N)^2}x^2 + 
        \frac{4}{C}\frac{\log \log N}{(\log N)^3}x^2 .
    \end{split}
\end{equation*}

 \begin{theorem}\label{T2}
Let $p>0$, \ $q_1>0$, \ $q_2 > 0$ be fixed with $p+ q_1+q_2 = 1$. 
Let $\mu(N)$ be the length of the longest at most $1+1$ contaminated run in $X_1,X_2,\ldots,X_N$.
Then for an integer $k$,
  \begin{equation}
           P(\mu(N)- [m(N)] < k) 
           =\exp\left(-p^{-(\log(C_0 p^{-2}q_1q_2) + H( k- \{m(N)\}) )}\right)
           \left(1+ \OO\left(\frac{1}{(\log N)^3}\right) \right) .          
   \end{equation}    
 \end{theorem}
%
\section{Proofs}  \label{proofs}
 Before proceeding with the proof, we shall consider fulfilment of some conditions given in the main Lemma 
 of \citep{csaki1987limit} for the case of $k=m$ (for fixed $m$) and $0<p<1$, and for some $\varepsilon>0$.
 \begin{remark}\label{rem1}
 First, we consider condition (SIII) and show that it is true for any large enough $m$.
 We have
\begin{equation}
    \begin{split}
        P(A_1) =& p^m + m(1-p)p^{m-1} + m(m-1)p^{m-2} q_1 q_2\\
        =& m(m-1)p^{m-2} q_1 q_2 \left \{ 1 + \frac{p(1-p)}{(m-1)q_1 q_2} + \frac{p^2}{m(m-1)q_1q_2}\right\} 
        \leq \frac{\varepsilon}{m} .
    \end{split}
\end{equation}
This inequality is true for any positive $\varepsilon$ if $m$ is large enough. \\
If  $m \approx \log N$, then $p^m \approx p^{\log N} = \frac{1}{N}$ and then, 
$\varepsilon \approx \frac{(\log N)^3}{N}$. (Here, $\log$ denotes logarithm to base $\frac{1}{p}$.)
 \end{remark}
 \begin{remark}\label{rem2}
     Now, consider condition (SII). If $i>m$, then $A_i$ and $A_1$ are independent, therefore
     \begin{equation}
         \sum_{i= m+1}^{2m} P(A_i|A_1) = m P(A_1) < \varepsilon ,
     \end{equation}
     which gives precisely the previous assumption in Remark  \ref{rem1}.
 \end{remark}
 \begin{lemma}\label{Lem1}
   Condition (SI) is satisfied  for $k=m$ in the following form 
   \begin{equation}
     |P(\bar{A_2}\Bar{A_3} \cdots \bar{A}_m| A_1) - \alpha| < \varepsilon  ,
   \end{equation}
     with 
     \begin{equation*}
     \alpha = \frac{C_0 +\frac{1}{m}C_1 +\frac{1}{m(m-1)}C_2}{1+ \frac{p(1-p)}{(m-1)q_1q_2} + \frac{p^2}{m(m-1)q_1q_2}} .
     \end{equation*}
\end{lemma}
\begin{proof}
     To begin, we shall be required to  divide the event $A_1$ into the following pairwise disjoint parts;
     \begin{equation*}
         A_1 = A_1^{(0)}\bigcup \left( \bigcup_{i=1}^m A_1^{(+)}(i)\right) 
         \bigcup \left( \bigcup_{i=1}^m A_1^{(-)}(i)\right) 
         \bigcup \left( \bigcup_{{i,j=1},{i\ne j}}^m A_1^{(2)}(i,j)\right),
     \end{equation*}
     where $A_1^{(0)}$ is the event that $X_1,X_2, \ldots, X_m$ is a pure run,\\
     $A_1^{(+)}(i)$ denotes that $X_i = +1$, while the rest are zeros,\\
     $A_1^{(-)}(i)$ denotes that $X_i = -1$, while the rest are zeros,\\
     finally, $A_1^{(2)}(i,j)$ denotes that $X_i = +1$, $X_j = -1$, while the rest are zeros. 
     
  We shall denote $X_i = +1$ by \coloredcircle{blue} while $X_j = -1$ by \coloredcircle{red} for ease of visualization.
     Then,
     \begin{equation*}
         \begin{split}
             P(A_1 \bar{A_2} \cdots \bar{A}_m) = & P(A_1^{(0)} \bar{A_2} \cdots\bar{A}_m) + 
             \sum_{i=1}^m P(A_1^{(+)}(i) \bar{A_2} \cdots\bar{A}_m) +\\
             & + \sum_{i=1}^m P(A_1^{(-)}(i) \bar{A_2} \cdots\bar{A}_m)+ 
             \sum_{i<j}^m P(A_1^{(2)}(i,j) \bar{A_2} \cdots\bar{A}_m) + \\
             & + \sum_{i>j}^m P(A_1^{(2)}(i,j) \bar{A_2} \cdots\bar{A}_m)\\
             := & Y^{(0)} + \sum_{i=1}^m  Y_i^{(+)} + \sum_{i=1}^m  Y_i^{(-)} 
             + \sum_{i<j}^m  Y_{i,j}^{(2)} + \sum_{i>j}^m  Y_{i,j}^{(2)}.
         \end{split}
     \end{equation*}
     Here, we can obtain the formula for $\sum_{i=1}^m  Y_i^{(-)}$ by interchanging the role of $q_1$ and $q_2$ 
     in the corresponding formula $\sum_{i=1}^m  Y_i^{(+)}$.\\
     Similarly, we can obtain $\sum_{i>j}^m  Y_{i,j}^{(2)}$ by interchanging the role of $q_1$ and $q_2$ 
     in the corresponding formula $\sum_{i<j}^m  Y_{i,j}^{(2)}$.
     
     Therefore,
     \begin{enumerate}[I.]
         \item  $Y^{(0)} = 0$, because the event is impossible.
         \item Now, calculate $ Y_i^{(+)} = P(A_1^{(+)}(i) \bar{A_2} \cdots\bar{A}_m) $. \\
     We want to evaluate probabilities corresponding to different values of $i$ as follows.
     
     \begin{enumerate}[(a)]
         \item If $i=1$, then the event is impossible. So $Y_1^{(+)} = 0$.
         \item  Let $1<i<m$, i.e\quad 
         $ \bigcirc,\cdots, \overset{i}{\coloredcircle{blue}}, \cdots,  \bigcirc, \cdots, \overset{m}{  \bigcirc}, 
          \overset{m+1}\square$.\\
  Then, the $m+1$ position should be $+1$. 
  Furthermore, on the positions $m+2,\ldots,m+i$,  it is not possible that all elements are zeros 
  and also not possible that there is a $-1$ and the rest of the elements are zeros. 
  So for this case,
  \begin{equation}\label{eqiib}
    Y_i^{(+)} = q_1 p^{m-1} q_1 \left (1 - p^{i-1} - (i-1) q_2 p^{i-2} \right ),  \quad \text{if} \quad 1< i<m.
  \end{equation}
  \item If $i=m$, i.e\quad 
  $ \bigcirc,\cdots, \bigcirc, \overset{m}{\coloredcircle{blue}}, \overset{m+1}{\coloredcircle{blue}}$,\\
  then $X_{m+1}$ should be $+1$ and other remaining elements to be arbitrary. 
  So this part is
  \begin{equation}\label{eqiic}
      Y_m^{(+)} = q_1 p^{m-1} q_1.
  \end{equation}
     \end{enumerate}
     \item Now, let us turn to $Y_{i,j}$, first we consider the case when $i<j$.
     \begin{enumerate}[(a)]
         \item  When $i=1$ and $j = m$, \quad 
         i.e \quad$\overset{1}{\coloredcircle{blue}}, \bigcirc,  \cdots,  \bigcirc,  
         \overset{m}{\coloredcircle{red}},\overset{m+1}\square$.\\
  Then, $X_{m+1}$ should be $-1$ and the remaining elements to be arbitrary. 
  So this part has
  \begin{equation}\label{eqiiia}
   Y_{1,m} = q_1 q_2 p^{m-2} q_2 .  
  \end{equation}
  \item Now, let $i=1$ and $j<m$, \quad i.e \quad
  $\overset{1}{\coloredcircle{blue}}, \bigcirc, \cdots, \overset{j}{\coloredcircle{red}},
  \cdots \overset{m}\bigcirc \overset{m+1}\square$.\\
  Then, $X_{m+1}$ should be $-1$. 
  Moreover, on positions $m+2, \ldots, m+j$, all elements being zeros is not possible, 
  neither is one +1 and the rest being zeros possible. 
  So
  \begin{equation}\label{eqiiib}
    Y_{1,j} = q_1 q_2 p^{m-2} q_2 \left (1 - p^{j-1} - p^{j-2}(j-1) q_1  \right ),  \quad \text{if} \quad 1<j<m. 
  \end{equation}
  \item Now, let $i>1$ and $j=m$ \quad i.e \quad 
  $\bigcirc,\cdots,\overset{i}{\coloredcircle{blue}},\bigcirc, \cdots, 
  \overset{m} {\coloredcircle{red}},\overset{m+1}\square$.\\
  Then, $X_{m+1}$ can either be a +1 or a $-1$. \\
  When $X_{m+1}$ is $-1$, then the remaining elements are arbitrary. 
  However, if $X_{m+1}$ is +1, then on positions $m+2, \ldots, m+i$, there should be at least one non-zero element. 
  So
  \begin{equation}\label{eqiiic}
    Y_{i,m} = q_1 q_2 p^{m-2} \left (q_1 \left (1 - p^{i-1}\right ) + q_2 \right ), \quad \text{if} \quad i>1 ,\quad j=m.  
  \end{equation}
  \item Consider now the case $i>1$ and $j<m$, \ i.e. \quad 
  $\bigcirc, \cdots, \overset{i}{\coloredcircle{blue}}, \bigcirc, 
  \cdots,\bigcirc, \overset{j}{\coloredcircle{red}}, \overset{m}{\bigcirc}$.\\
  We divide this event into two parts.\\
  First, let $X_{m+1}= +1$, \quad 
  $\bigcirc, \cdots, \overset{i}{\coloredcircle{blue}}, \cdots,\bigcirc, \overset{j}{\coloredcircle{red}}, 
  \dots,\overset{m}{\bigcirc}, \overset{m+1}{\coloredcircle{blue}}$.\\
  Then, it is not possible that the elements in positions $m+2, \ldots, m+i$ are all zeros. 
  Neither is it also possible that there is one $-1$ among $m+2, \ldots, m+i$ while all $m+i+1, \ldots, m+j$ are zeros.
  So this part of $Y_{i,j}$ is
  \begin{equation}\label{eqiiid1}
  q_1 q_2 p^{m-2} q_1 \left (1 - p^{i-1} - (i-1)  p^{i-2} q_2  p^{j-i}\right ),  
  \quad \text{if} \quad i>1,\quad j<m.   
  \end{equation}
   Finally, now let $X_{m+1}= -1$. \quad \quad 
   $\bigcirc, \cdots, \overset{i}{\coloredcircle{blue}}, \cdots,\bigcirc, \overset{j}{\coloredcircle{red}}, 
   \dots,\overset{m}{\bigcirc}, \overset{m+1}{\coloredcircle{red}}.$\\
  Then it is not possible that all elements in $m+2, \ldots, m+j$ are zeros, 
  neither is it possible that among $m+2, \ldots, m+j$ there is one +1 and the rest are zeros. 
  So this second part of $Y_{i,j}$ is 
  \begin{equation}\label{eqiiid2}
   q_1 q_2 p^{m-2} q_2 \left (1 - p^{j-1} - (j-1) q_1 p^{j-2} \right ),  \quad \text{if} \quad i>1,\quad j<m.   
  \end{equation}
     \end{enumerate}
     \end{enumerate}
    
  Summing equations \eqref{eqiib} and \eqref{eqiic}, we get $Y_i^{(+)}$ and 
  consequently by interchanging the roles of $q_1$ and $q_2$ we obtain $Y_j^{(-)}$ as follows
  \begin{equation}\label{eqbc}
      \begin{split}
      \sum_{i=1}^m  Y_i^{(+)} + \sum_{i=1}^m  Y_i^{(-)} =& 
      \sum_{i=2}^{m-1} q_1^2 p^{m-1} \left (1 - p^{i-1} - (i-1) q_2 p^{i-2} \right ) + q_1^2 p^{m-1}  +\\
      & +\sum_{i=2}^{m-1} q_2^2 p^{m-1} \left (1 - p^{i-1} - (i-1) q_1 p^{i-2} \right ) + q_2^2 p^{m-1}\\
      = & (m-1)p^{m-1}\left( q_1^2 + q_2^2\right) - \left( q_1^2 + q_2^2\right)p^{m-1}\sum_{i=2}^{m-1} p^{i-1} - \\
      & - p^{m-1}\left( q_1^2 q_2 + q_2^2 q_1\right)\sum_{i=2}^{m-1} (i-1)p^{i-2}\\
      = & (m-1)p^{m-1}\left( q_1^2 + q_2^2\right) - \left( q_1^2 + q_2^2\right)p^{m-1} \frac{p^{m-1}- p}{p-1} - \\
      - & p^{m-1} q_1 q_2 \left( q_1 + q_2\right) \left( \frac{(m-2) p^{m-2} -1}{p-1} + \frac{p-p^{m-2}}{(p-1)^2}\right) =\\
      = & p^{m-1} \left\{ (q_1^2 +q_2^2) \left [(m-1) - \frac{p^{m-1}}{p-1} + \frac{p}{p-1}\right ] \right.+\\
      & + \left. q_1 q_2 \left[ (m-2) p^{m-2} -1 + \frac{p}{p-1} -\frac{p^{m-2}}{p-1}\right] \right\} =\\
      = & p^{m-1} \left \{ (q_1^2 + q_2^2) \left[ (m-1) + \frac{p}{p-1}\right] + q_1 q_2 \frac{1}{p-1} \right. -\\
      & - \left. (q_1^2 + q_2^2) 
      \frac{p^{m-1}}{p-1} + q_1 q_2 \left[ (m-2) p^{m-2} - \frac{p^{m-2}}{p-1}\right]\right \} .
      \end{split}
  \end{equation}
  Here above, we applied
  \begin{equation*}
      \sum_{i=a}^b i p^{i-1} = \frac{bp^b-ap^{a-1}}{p-1} + \frac{p^a-p^b}{(p-1)^2} ,
  \end{equation*}
    which can be obtained by differentiating the known formula for the sum of a geometric sequence.
  Similarly, summing equations \eqref{eqiiia}, \eqref{eqiiib}, \eqref{eqiiic}, \eqref{eqiiid1} and \eqref{eqiiid2}
   together with their corresponding changed versions got by interchanging roles of $q_1$ and $q_2$, we obtain  
  \begin{equation*}
      \begin{split}
          \sum_{i<j} & Y_{i,j} + \sum_{i>j} Y_{i,j} =  q_1q_2p^{m-2}(q_1+q_2) +\\
          & + q_1q_2p^{m-2} \sum_{j=2}^{m-1} \left( q_2 \left( 1-p^{j-1}-p^{j-2} (j-1)q_1\right) + q_1 \left( 1-p^{j-1} -p^{j-2} (j-1)q_2\right)\right) +\\
          & + q_1 q_2 p^{m-2} \left[ \sum_{i=2}^{m-1} \left( q_1 \left( 1-p^{i-1}\right) + q_2\right) + \sum_{i=2}^{m-1} \left( q_2 \left( 1-p^{i-1}\right) + q_1\right)\right] +\\
          & + \sum_{i=2}^{m-2}\sum_{j=i+1}^{m-1}\left[q_1 q_2 p^{m-2} q_1 \left( 1-p^{i-1} - (i-1)p^{i-2} q_2 p^{j-i}\right) \right.+ \\
          & \left.+ q_1q_2 p^{m-2} q_2 \left( 1-p^{j-1}- (j-1)q_1 p^{j-2}\right)
          + q_1q_2 p^{m-2} q_2 \left( 1-p^{i-1}- (i-1)p^{i-2} q_1 p^{j-i}\right)\right.+ \\
          & \left.+ q_1q_2 p^{m-2} q_1 \left( 1-p^{j-1}- (j-1) q_2 p^{j-2} \right)\right] =\\
          & = q_1q_2 p^{m-2} \left\{ (q_1 +q_2) \left( m-1 - \frac{p^{m-1}-p}{p-1}\right) - 2q_1q_2 \sum_{j=2}^{m-1} (j-1) p^{j-2} \right.+ \\
          & \left.+ 2(q_1 +q_2)(m-2)- (q_1 +q_2)\frac{p^{m-1}-p}{p-1} + \sum_{i=2}^{m-2}\sum_{j=i+1}^{m-1} \left[ (q_1 +q_2) (1-p^{i-1}) \right.\right.-\\
          & \left.\left. - 2q_1 q_2 (i-1) p^{j-2} + (q_1 +q_2) (1-p^{j-1})-2q_1 q_2 (j-1) p^{j-2} \right]\right\} =\\
          & = q_1q_2 p^{m-2} \left\{ (q_1 +q_2)(3m-5) -2(q_1 +q_2) \frac{p^{m-1}-p}{p-1}\right.- \\
          & \left. -2 q_1 q_2 \left( \frac{(m-2)p^{m-2}-1}{p-1} + \frac{p-p^{m-2}}{(p-1)^2}\right) \right. +\\
          & \left.+ \sum_{i=2}^{m-2} \left[ (q_1 +q_2)(1-p^{i-1})(m-i-1) - 2q_1 q_2 (i-1) \frac{p^{m-2}-p^{i-1}}{p-1} \right.\right. +\\
          & \left.\left.+ (q_1 +q_2)(m-i-1) -(q_1 +q_2) \frac{p^{m-1}-p^i}{p-1} \right.\right.-\\
          & \left.\left.-2q_1 q_2 \left( \frac{(m-2) p^{m-2}- i p^{i-1}}{p-1} + \frac{p^i-p^{m-2}}{(p-1)^2}\right)\right] \right\}   .      
      \end{split}
  \end{equation*}
  \begin{equation*}
      \begin{split}
          \sum_{i<j} & Y_{i,j} + \sum_{i>j} Y_{i,j} ={}\\
          &=q_1q_2 p^{m-2}\left\{ (q_1 +q_2)(3m-5) + 2(p^{m-1}-p) -2q_1 q_2 \left( \frac{(m-2) p^{m-2}- 1}{p-1} + 
          \frac{p-p^{m-2}}{(p-1)^2}\right)\right.+\\
          & \left. + (q_1 +q_2)(m-3) \frac{m-2}{2} - (q_1 +q_2)m \sum_{i=2}^{m-2}p^{i-1} + (q_1 +q_2)\sum_{i=2}^{m-2}p^{i-1}(i+1) \right.-\\
          & \left.- 2q_1 q_2 \frac{p^{m-2}}{p-1} (m-3) \frac{m-2}{2} + 2\frac{q_1 q_2}{p-1}\sum_{i=2}^{m-2} (i-1)p^{i-1} + (q_1 +q_2)(m-3)\frac{m-2}{2} \right.-\\
          & \left.- (q_1 +q_2)\frac{p^{m-1}}{p-1} (m-3) + \frac{q_1 +q_2}{p-1}\sum_{i=2}^{m-2}p^i- 
          \frac{2q_1q_2}{p-1} (m-2)p^{m-2}(m-3) + \frac{2q_1q_2}{p-1}\sum_{i=2}^{m-2} i p^{i-1}\right.-\\
          & \left.-\frac{2q_1q_2}{(p-1)^2}\sum_{i=2}^{m-2}p^i + \frac{2q_1q_2}{(p-1)^2} p^{m-2}(m-3)\right\}=\\
          & = 
          q_1q_2 p^{m-2}\left\{ (q_1 +q_2)(3m-5 +(m-3)(m-2)) + 2\left(p^{m-1}-p\right) \right.-\\
          & \left.- 2q_1q_2 \frac{(m-2)p^{m-2}}{p-1} + 2q_1q_2 \frac{1}{p-1} - \frac{2q_1q_2 p}{(p-1)^2} + 
          \frac{2q_1q_2 p^{m-2}}{(p-1)^2} -(q_1+q_2)m \frac{p^{m-2}-p}{p-1} \right.+\\
          & \left.+ (q_1+q_2) \left( \frac{1}{p} .\frac{(m-1)p^{m-1}-3p^2}{p-1} + \frac{1}{p}
           \frac{p^3-p^{m-1}}{(p-1)^2}\right) -2q_1q_2 \frac{p^{m-2}}{p-1} \frac{(m-3)(m-2)}{2} \right.+\\
          & \left.+ \frac{2q_1q_2}{p-1} \left( \frac{(m-3)p^{m-3}-1}{p-1} + \frac{p-p^{m-3}}{(p-1)^2}\right)p + 
          p^{m-1}(m-3)-\frac{p^{m-1}-p^2}{p-1} \right.-\\
          & \left.- \frac{ 2q_1q_2 p^{m-2}}{p-1}(m-2)(m-3) + \frac{2q_1q_2}{p-1} \left( \frac{(m-2)p^{m-2}-2p}{p-1} 
          + \frac{p^2-p^{m-2}}{(p-1)^2}\right) \right.-\\
          & \left.- \frac{2q_1q_2}{(p-1)^2} \frac{p^{m-1}-p^2}{p-1} + \frac{2q_1q_2}{(p-1)^2} p^{m-2}(m-3)\right\}=\\
          & = 
          q_1q_2 p^{m-2}\left\{ (q_1 +q_2)(3m-5 + m^2-5m+6) + 2p^{m-1}-2p-2q_1q_2 \frac{(m-2)p^{m-2}}{p-1} \right.+\\
          & \left.+2q_1q_2 \frac{1}{p-1} - 2q_1q_2 \frac{p}{(p-1)^2}+ 2q_1q_2 \frac{p^{m-2}}{(p-1)^2} + mp^{m-2}- pm-(m-1)p^{m-2}  \right.+\\
          & \left.+3p -\frac{p^2}{p-1}+\frac{p^{m-2}}{p-1} -q_1q_2 \frac{p^{m-2}}{p-1}(m-2)(m-3)+ \frac{2q_1q_2(m-3)p^{m-2}}{(p-1)^2} - \frac{2q_1q_2 p}{(p-1)^2} \right.+\\
          & \left.+ \frac{2q_1q_2 p^2}{(p-1)^3}- \frac{2q_1q_2 p^{m-2}}{(p-1)^3}+ p^{m-1}(m-3) - \frac{p^{m-1}}{p-1}+ \frac{p^2}{p-1} -\frac{2q_1q_2 p^{m-2}}{p-1}(m-2)(m-3)\right.+\\
          & \left.+\frac{2q_1q_2 (m-2)p^{m-2}}{(p-1)^2} -\frac{4q_1q_2 p}{(p-1)^2} + \frac{2q_1q_2 p^2}{(p-1)^3} - \frac{2q_1q_2 p^{m-2}}{(p-1)^3}- \frac{2q_1q_2 p^{m-1}}{(p-1)^3}\right.+\\
          & \left. +  \frac{2q_1q_2 p^2}{(p-1)^3} + \frac{2q_1q_2 p^{m-2}(m-3)}{(p-1)^2}\right\} .
      \end{split}
  \end{equation*}

  \begin{equation}\begin{split}\label{eqiii}
    \sum_{i<j}& Y_{i,j} + \sum_{i>j} Y_{i,j} = {}\\
    & q_1q_2 p^{m-2}\left\{ (q_1 +q_2)(m-1)^2 - mp+ p+ 2q_1q_2 \frac{1}{p-1}- 2q_1q_2 \frac{p}{(p-1)^2}- 
    \frac{2q_1q_2 p}{(p-1)^2}\right.+\\
          & \left.+ \frac{2q_1q_2 p^2}{(p-1)^3}+ \frac{p^2}{p-1}- \frac{4q_1q_2 p}{(p-1)^2}+ \frac{2q_1q_2 p^2}{(p-1)^3}
           + \frac{2q_1q_2 p^2}{(p-1)^3}- \frac{p^2}{p-1} - 2q_1q_2\frac{(m-2) p^{m-2}}{p-1}\right.+\\
          & \left.+ mp^{m-2}-(m-1)p^{m-2}-q_1q_2\frac{p^{m-2}}{p-1}(m-3)(m-2)+
          \frac{2q_1q_2 (m-3)p^{m-2}}{(p-1)^2}\right.+\\
          & \left.+ p^{m-1}(m-3) - \frac{2q_1q_2p^{m-2}}{p-1}(m-2)(m-3) +\frac{2q_1q_2(m-2)p^{m-2}}{(p-1)^2} 
          +\frac{2q_1q_2p^{m-2}(m-3)}{(p-1)^2}\right.+\\
          & \left.+ 2 p^{m-1} + \frac{2q_1q_2p^{m-2}}{(p-1)^2} + \frac{p^{m-2}}{p-1} -  \frac{2q_1q_2p^{m-2}}{(p-1)^3} 
          -\frac{p^{m-1}}{p-1}-  \frac{2q_1q_2p^{m-2}}{(p-1)^3} -\frac{2q_1q_2p^{m-1}}{(p-1)^3}\right\} =\\
   &= 
   q_1q_2 p^{m-2}\left\{ (q_1 +q_2)(m-1)^2-p(m-1)+ \frac{2q_1q_2}{(p-1)^3} \left[ (p^2-2p+1)-4(p^2-p)\right.\right.+\\
          & \left.\left. + 3p^2\right] - 3q_1q_2 \frac{p^{m-2}}{p-1}(m-3)(m-2) + p^{m-1}(m-3) + 
          \frac{2q_1q_2}{p-1}p^{m-3}\left( -p(p-1)(m-2) \right.\right.+\\
          & \left.\left.+ p(m-3) + (m-2)p + p(m-3)\right) + p^{m-2} + 2p^{m-1} + \frac{p^{m-2}}{p-1}- 
          \frac{p^{m-1}}{p-1}\right.+\\
          & \left. +\frac{2q_1 q_2 p^{m-3}}{(p-1)^3} [p(p-1)-p-p-p^2] \right\} = \\         
    & = 
    q_1q_2 p^{m-2} \left\{ (q_1 +q_2)(m-1)^2-p(m-1)+ \frac{2q_1q_2}{(p-1)^3}(2p+1)\right.-\\
          & \left.- 3q_1q_2\frac{p^{m-2}}{p-1}(m-2)(m-3) + p^{m-1}(m-3) + \frac{2q_1q_2}{(p-1)^2}p^{m-3} p 
          \left[(4-p)m+2p-10 \right] \right.+\\
          & \left.+ \frac{p^{m-2}}{p-1}(2p(p-1)) + \frac{2q_1q_2}{(p-1)^3}p^{m-3}(-3p)\right\}=\\   
      &=
      q_1q_2 p^{m-2} \left\{ (q_1 +q_2)(m-1)^2-p(m-1) + 2(2p+1)\frac{q_1q_2}{(p-1)^3}\right.-\\
          & \left. -3q_1q_2\frac{p^{m-2}}{p-1}(m-2)(m-3) + p^{m-1}(m-3) + \frac{2q_1q_2}{(p-1)^2}p^{m-3} p 
          \left[(4-p)m+2p-10 \right]\right.+\\
          & \left.+ 2p^{m-1} - \frac{6q_1q_2p^{m-2}}{(p-1)^3}\right\}=\\
          &= 
          q_1q_2 p^{m-2}\left\{ \underbrace{m(m-1)(q_1+q_2)-(m-1)}_{m^2(q_1+q_2)-m(q_1+q_2+1)+1} + 
          \frac{2(2p+1)q_1q_2}{(p-1)^3} \right.+\\
          & \left.+ \frac{q_1q_2p^{m-2}}{(p-1)^3} \Big( -3(p-1)^2 m^2 + m(p-1)(13p-7) + 
          (-14p^2+12p-4)\Big) +p^{m-1}(m-1)\right\} .  
  \end{split}
  \end{equation}
  Therefore, combining \eqref{eqbc} and \eqref{eqiii} and by some simplification, we obtain
  \begin{equation*}
      \begin{split}
          P(A_1 \bar{A_2} \cdots\bar{A}_m) {}& = p^{m-1}\left\{(q_1^2+q_2^2)\left[ (m-1)+ 
          \frac{p}{p-1}\right] +q_1q_2\frac{1}{p-1}+ \OO(mp^m)\right\}+\\
          & + q_1q_2p^{m-2}\left\{ m(m-1)(q_1+q_2)- (m-1)+ \frac{2(2p+1)q_1q_2}{(p-1)^3} + \OO(m^2p^m)\right\}=\\
          = {}& m(m-1)p^{m-2}q_1q_2\left\{ \frac{p(q_1^2+q_2^2)}{mq_1q_2} + 
          \frac{(q_1^2+q_2^2)p^2}{m(m-1)q_1q_2(p-1)} \right.+\\
          & \left.+ \frac{p}{m(m-1)(p-1)} + \OO\left(\frac{p^m}{m}\right) + (q_1+q_2) -\frac{1}{m} \right.+\\
          & \left.+ \frac{2(2p+1)q_1q_2}{m(m-1)(p-1)^3} + \OO(p^m)\right\}=\\
          = {}& m(m-1)p^{m-2}q_1q_2\left\{ C_0 + \frac{1}{m}C_1 + \frac{1}{m(m-1)}C_2 + \OO(p^m)\right\},
      \end{split}
  \end{equation*}
  where $C_0=(q_1+q_2)$, $C_1= \frac{p(q_1^2+q_2^2)}{q_1q_2}-1$ , 
  $C_2= \frac{(q_1^2+q_2^2)p^2}{q_1q_2(p-1)} + \frac{p}{p-1} + \frac{2(2p+1)q_1q_2}{(p-1)^3}$. \\
  So,
  \begin{equation*}
   P(\bar{A_2} \cdots\bar{A}_m|A_1) = \frac{C_0 + \frac{1}{m}C_1 +\frac{1}{m(m-1)}C_2 + \OO(p^m) }{1+ 
   \frac{p(1-p)}{(m-1)q_1q_2} + \frac{p^2}{m(m-1)q_1q_2}} .
  \end{equation*}
  We therefore satisfy Lemma \ref{Lem1}.
  \begin{equation*}
      |P(\bar{A_2} \cdots\bar{A}_m|A_1)-\alpha|< \varepsilon ,
  \end{equation*}
  where
  \begin{equation*}
      \alpha = \frac{C_0 +\frac{1}{m}C_1 +\frac{1}{m(m-1)}C_2  }{1+ \frac{p(1-p)}{(m-1)q_1q_2} + 
      \frac{p^2}{m(m-1)q_1q_2}} 
  \end{equation*}
  and $\varepsilon= \OO(p^m)$.
  \end{proof}
 %
 \begin{proof}[\textbf{Proof of Theorem \ref{T1}}]
 We shall apply the Main Lemma (stationary case, finite form) in \citep{csaki1987limit} since its conditions in 
 Remarks \ref{rem1}, \ref{rem2} and Lemma \ref{Lem1} are satisfied. 
 
 As $\tau_m$ is the first hitting time of the at most $1+1$ contaminated run of length $m$, we have
 \begin{equation*}
     P(\tau_m \alpha P(A_1)> x)= P\left( \tau_m > \frac{x}{\alpha P(A_1)}\right)= P(\tau_m>N) ,
 \end{equation*}
 where $N$ is the integer part of $\frac{x}{\alpha P(A_1)}$.\\
 We can see that 
 $ P(\tau_m>N) = P(\bar{A_1} \cdots\bar{A}_{N_1})$,  where in this case, $N_1= N-m+1$.\\
 By \citep{csaki1987limit} main lemma,
 \begin{equation*}
   \begin{split}
       e^{-(\alpha+10\varepsilon)N_1P(A_1)-2mP(A_1)} <P(\bar{A_1}, \cdots,\bar{A}_{N_1}) < 
       e^{-(\alpha-10\varepsilon)N_1P(A_1)+2mP(A_1)} .
   \end{split}  
 \end{equation*}
 Let $\varepsilon_0>0$ be fixed and choose $\varepsilon$ so that $10\varepsilon=\frac{\varepsilon_0}{m}$. 
 We can do this because $\varepsilon=O(p^m)$.\\
 Now
 \begin{equation*}
     \begin{split}
      P(\tau_m>N)=  &P(\bar{A_1} \cdots\bar{A}_{N_1}) \sim  
      e^{-(\alpha\pm \frac{\varepsilon_0}{m})N_1P(A_1)\pm2mP(A_1)}\sim  \\
      & \sim 
   e^{-\left(\alpha\pm \frac{\varepsilon_0}{m}\right)\left( \frac{x}{\alpha P(A_1)}-m+1\right) P(A_1)}e^{\pm2mP(A_1)}=\\
      &= e^{-\left(\alpha\pm \frac{\varepsilon_0}{m}\right)\left( \frac{xP(A_1)}{\alpha P(A_1)}\right)}
       e^{-\left(\alpha\pm \frac{\varepsilon_0}{m}\right)(-m+1)P(A_1)} e^{\pm2mP(A_1)}\sim\\
      &\sim e^{-x}
     \end{split}
 \end{equation*}
 because $\alpha\approx C_0$, where $C_0=(q_1+q_2)$, $mP(A_1)\to 0$ as $m\to \infty$.\\
 So we obtained the statement of Theorem \ref{T1}
 \begin{equation*}
   P(\tau_m \alpha P(A_1)> x)\sim e^{-x} \quad {\text {as}} \quad m\to \infty.  
 \end{equation*}
 \end{proof}
 %
 \begin{proof}[\textbf{Proof of Theorem \ref{T2}}]
 Let $N_1=N-m+1$, where $m$ will be specified so that $m\sim \log N$.
Then,
\begin{equation*}
    \begin{split}
        P(\mu(N) < m) & =  P(\bar{A_1} \cdots \bar{A}_{N_1}) \sim 
        e^{-(\alpha\pm 10\varepsilon)N_1P(A_1) \pm 2m P(A_1)} \\
        & = e^{-\alpha N_1 P(A_1)} e^{\pm 10 \varepsilon N_1 P(A_1)} e^{\pm 2m P(A_1)} .
    \end{split}
\end{equation*}
As $mP(A_1) \sim m^3 p^m \sim \frac{(\log N)^3}{N}$,
 so $e^{\pm 2m P(A_1)} = 1 + \OO\left(\frac{(\log N)^3}{N} \right)$.
Similarly, as $\varepsilon = \OO\left(p^m\right)$ and $m \approx \log N$,
\begin{equation*}
    e^{\pm 10 \varepsilon N_1 P(A_1)} \sim e^{\pm 10 P(A_1)} 
    \sim e^{\pm c(\log N)^2/N} = 1 + O \left( \frac{(\log N)^2}{N}\right) .
\end{equation*}
Therefore, we can calculate
\begin{equation*}
    e^{-\alpha N_1 P(A_1)} = e^{-\alpha N P(A_1)} 
    \underbrace{e^{+ \alpha (m-1)P(A_1)}}_{1 +    \OO \left( \frac{(\log N)^3}{N}\right)} .
\end{equation*}
So we have to calculate:
     $$
     e^{-\alpha N P(A_1)} = e^{-l}, 
     $$
 where
\begin{equation*}
\begin{split}
    l = {}&\alpha N P(A_1) \\
= {}& \frac{C_0 + \frac{1}{m} C_1 + \frac{1}{m(m-1)} C_2}{1 + \frac{p(1-p)}{(m-1)q_1 q_2} + \frac{p^2}{m(m-1)q_1 q_2}}
 N m(m-1)p^{m-2}(q_1 q_2) \left( 1 + \frac{p(1-p)}{(m-1)q_1 q_2} + \frac{p^2}{m(m-1)q_1 q_2}\right)\\
    & = N p^{m-2} q_1 q_2 \left( m(m-1)C_0 + (m-1)C_1 + C_2\right)\\
    & = N p^{m-2} q_1 q_2 \left( m^2C_0 + m(C_1-C_0) + C_2 - C_1\right).
\end{split}
\end{equation*}
 
 Our aim is to find $m(N)$ so that the asymptotic behaviour of $P(\mu(N)-[m(N)]<k)$ can be obtained. 
 Here $[m(N)]$ is the integer part of $m(N)$ and $\{m(N)\}$ denotes its fractional part. 
 Then
 $$
 P(\mu(N)-[m(N)]<k)=P(\mu(N)<m(N)+k-\left\{m(N) \right\}).
 $$
 Let $ m= m(N)+k-\left\{m(N) \right\}$. 
 So
 $$
 P(\mu(N)-[m(N)]<k)=P(\mu(N)<m)= e^{-l}\left( 1+ \OO\left( \frac{(\log N)^3}{N}\right)\right),
 $$
 where $\log $ is the logarithm to base $1/p$.
 We want to find $m(N)$ so that the remainder term in the exponent $l$ be small.
 We shall do it step by step using several Taylor's expansions.
 We try to find $m(N)$ as $\log N+A$, where $A$ will be specified later. 
 So $m= \log N+A+k-\{ m(N)\}$.
 
 Then, using Taylor expansion 
 $\log (x_0+y)= \log x_0 + \frac{y}{C x_0}-\frac{1}{2C}\frac{y^2}{x_0^2} + \frac{1}{3C}\frac{y^3}{\tilde{x}^3}$, 
 where $\tilde{x}$ is between $x_0$ and $x_0+y$ and $x_0>0$, $x_0+y>0$, and
 where $C=\ln\frac{1}{p}$ and $\ln $ is the logarithm to base $e$, we obtain
 \begin{equation*}
     \begin{split}
         L= {}&\log l= \log \left( p^{-2}q_1q_2\right)-m + \log N + \log (m^2 C_0) + 
         \frac{m(C_1-C_0)+ (C_2-C_1)}{C m^2C_0} -\\
         -{}& \frac{1}{2C} \frac{(m(C_1-C_0) + (C_2-C_1))^2}{(C m^2 C_0)^2} + \OO\left(\frac{1}{m^3}\right)\\
         ={}& \log \left( p^{-2}q_1q_2\right)- m + \log N + 2 \log m + \log C_0 + \frac{C_1-C_0}{C C_0 m}+ 
         \frac{1}{m^2} K + \OO\left( \frac{1}{m^3}\right),
 \end{split}
 \end{equation*}        
     where \quad
     $ K = \frac{C_2-C_1}{CC_0} - \frac{(C_1-C_0)^2}{2CC_0^2} = \frac{2C_0C_2 - C_1^2 - C_0^2}{2CC_0^2}$.
  Now as \ $m = \log N + A +k- \left\{ m(N)\right\}$, \ so 
  \begin{equation*}
     \begin{split}
        L ={}& \log\left( C_0 p^{-2}q_1q_2\right)- \log N - A -\left( k- \{m(N)\}\right) + \log N + \\ 
        & + 2\log \left( \log N + A + k- \{m(N)\}\right) 
         + \frac{C_1-C_0}{C C_0 m} + \frac{1}{m^2}K + \OO\left( \frac{1}{m^3}\right) .
     \end{split}
 \end{equation*}
Again applying the Taylor expansion of 
$\log(x_0+y)= \log x_0 + \frac{y}{Cx_0}-\frac{1}{2C}\frac{y^2}{ x_0^2} + \frac{1}{3C}\frac{y^3}{\tilde{x_0}^3}$, 
we have
 \begin{equation*}
     \begin{split}
         L= {}&\log(C_0 p^{-2} q_1 q_2) - A -(k-\left\{ m(N)\right\}) + 2 \left(\log \log N + 
         \frac{A +k -\left\{ m(N)\right\} }{C\log N} - \right.\\
         & \left.-\frac{1}{2C} \frac{(A +k -\left\{ m(N)\right\})^2}{(\log N)^2} + 
         \frac{1}{3C} \frac{(A +k -\left\{ m(N)\right\})^3}{(\log N)^3}\right) +\\
         & + \OO\left( \frac{1}{(\log N)^3}\right) + \frac{C_1-C_0}{C C_0 m} + \frac{K}{m^2} 
         + \OO\left( \frac{1}{m^3}\right)\\
         = {}& \log (C_0 p^{-2} q_1 q_2) - A -(K-\left\{ m(N)\right\}) + 2 \log \log N + 
         \frac{2(A +k -\left\{ m(N)\right\}) }{C\log N}-\\
         - {}&\frac{1}{C} \frac{(A +k -\left\{ m(N)\right\})^2 }{(\log N)^2} + 
         \frac{2}{3C} \frac{(A +k -\left\{ m(N)\right\})^3 }{(\log N)^3} + \frac{C_1-C_0}{C C_0 m} + 
         \frac{K}{m^2} + \OO\left( \frac{1}{(\log N)^3}\right).
     \end{split}
 \end{equation*}
 Now, we let $A = 2 \log \log N + B$. Then,
   \begin{equation*}
       \begin{split}
        L=  {}&\log(C_0 p^{-2} q_1 q_2)- 2 \log \log N - B - (k-\left\{ m(N)\right\}) + 2 \log \log N +\\
        + {} &\frac{4 \log \log N + 2 B + 2(k-\left\{ m(N)\right\})}{C \log N} - 
        \frac{1}{C} \frac{A^2 + 2A(k-\left\{ m(N)\right\}) + (k-\left\{ m(N)\right\})^2}{(\log N)^2}+\\
        + {} & \frac{2}{3C} \left(\frac{A^3}{(\log N)^3} + \frac{3 A^2(k-\left\{ m(N)\right\})}{(\log N)^3} + 
        \frac{3 A(k-\left\{ m(N)\right\})^2}{(\log N)^3} + \frac{(k-\left\{ m(N)\right\})^3}{(\log N)^3}\right) +\\
        + {} & \frac{C_1-C_0}{C C_0 m} + \frac{K}{m^2} + \OO\left( \frac{1}{(\log N)^3}\right)\\
        ={}& \log(C_0 p^{-2} q_1 q_2) - B- (k-\left\{ m(N)\right\}) + \frac{4 \log \log N}{C \log N} + 
        \frac{2B}{C \log N} + \frac{2 (k-\left\{ m(N)\right\})}{C \log N}-\\
        -{}& \frac{1}{C} \frac{A^2}{(\log N)^2} - \frac{2}{C}\frac{A(k-\left\{ m(N)\right\})}{(\log N)^2} - 
        \frac{1}{C}\frac{(k-\left\{ m(N)\right\})^2}{(\log N)^2} + \frac{2}{3C} \frac{A^3}{(\log N)^3} +\\
        +{}& \frac{2A^2(k-\left\{ m(N)\right\})}{C(\log N)^3} + \frac{2A(k-\left\{ m(N)\right\})^2}{C(\log N)^3} +
         \frac{C_1-C_0}{C C_0 m} + \frac{K}{m^2} + \OO\left( \frac{1}{(\log N)^3}\right).
       \end{split}
   \end{equation*}
Let $B = \frac{4 \log \log N}{C \log N} + D$. 
Then,
   \begin{equation*}
       \begin{split}
        L=  {}&\log(C_0 p^{-2} q_1 q_2)- D -\left( k- \{m(N)\}\right) + \frac{8\log\log N}{C^2(\log N)^2} + 
        \frac{2D}{C\log N} + \frac{2(k- \{m(N)\})}{C\log N} -\\
        -{}& \frac{1}{C}\frac{(2\log\log N + B)^2}{(\log N)^2}- \frac{2}{C}\frac{(2\log\log N + B)
        (k-\left\{ m(N)\right\})}{(\log N)^2} - \frac{1}{C}\frac{(k-\left\{ m(N)\right\})^2}{(\log N)^2}+ \\
        + {}& \frac{2}{3C} \frac{(2\log\log N + B)^3}{(\log N)^3} + 
        \frac{2(2\log\log N + B)^2(k- \{m(N)\})}{C(\log N)^3} +\\
        +{}& \frac{2(2\log\log N + B)(k- \{m(N)\})^2}{C(\log N)^3} +\frac{C_1-C_0}{C C_0 m} + 
        \frac{K}{m^2} +  \OO\left( \frac{1}{(\log N)^3}\right).
       \end{split}
   \end{equation*}
   We shall use the Taylor expansion of the function $\frac{1}{x}$ as
   \begin{equation*}
       \frac{1}{x_0+x}= \frac{1}{x_0} - \frac{x}{x_0^2} + \frac{x^2}{x_0^3} - \frac{x^3}{\tilde{x}^4} ,
   \end{equation*}
    where $\tilde{x}$ is between $x_0$ and $x_0+x$ and where $x_0>0$ and $x_0 + x >0$.
   Since $m = \log N + A + k - \left\{ m(N)\right\}$, so
   \begin{equation*}
     \frac{1}{m} = \frac{1}{\log N} - \frac{A +k- \{m(N)\}}{(\log N)^2} + 
     \frac{(A +k- \{m(N)\})^2}{(\log N)^3}+ \OO\left( \frac{1}{(\log N)^3}\right).  
   \end{equation*}
   and
   \begin{equation*}
   \begin{split}
     \frac{1}{m^2} = {}&\frac{1}{(\log N)^2 + 2\log N(A +k- \{m(N)\}) + (A +k- \{m(N)\})^2 } =\\
     ={}& \frac{1}{(\log N)^2} - \frac{2 \log N (A +k- \{m(N)\})}{(\log N)^4} + \OO\left( \frac{1}{(\log N)^3}\right).
   \end{split}      
   \end{equation*}
   Now let 
   \begin{equation*}
       D = \frac{C_1-C_0}{C C_0} \frac{1}{\log N} + E.
   \end{equation*}
   Then,
   \begin{equation*}
   \begin{split}  
      L=  {}&\log(C_0 p^{-2} q_1 q_2) - \frac{C_1-C_0}{C C_0} \frac{1}{\log N} - E-\left( k- \{m(N)\}\right) +
       \frac{8\log\log N}{C^2(\log N)^2} +  \\
      +{} & \frac{2(C_1-C_0)}{C^2 C_0} \frac{1}{(\log N)^2} +\frac{2E}{C \log N} + 
      \frac{2( k- \{m(N)\}}{C \log N} - \frac{4}{C} \frac{(\log\log N)^2}{(\log N)^2} - 
      \frac{4B}{C} \frac{\log\log N}{(\log N)^2} -  \\
      - {}& \frac{4}{C} \frac{(\log\log N)(k-\left\{ m(N)\right\})}{(\log N)^2} -
      \frac{2B}{C} \frac{(k-\left\{ m(N)\right\})}{(\log N)^2} -  
      \frac{1}{C} \frac{(k-\left\{ m(N)\right\})^2}{(\log N)^2} +  
      \frac{16}{3C} \frac{(\log \log N)^3}{(\log N)^3} + \\
      + {}& \frac{8B}{C} \frac{(\log \log N)^2}{(\log N)^3} +
      \frac{8}{C} \frac{(\log\log N)^2(k-\left\{ m(N)\right\})}{(\log N)^3} + 
      \frac{4}{C} \frac{(\log\log N)(k-\left\{ m(N)\right\})^2}{(\log N)^3} +\\
      + {}& \frac{C_1-C_0}{C C_0}\left( \frac{1}{\log N} - \frac{A +k- \{m(N)\}}{(\log N)^2} + 
      \frac{(2\log \log N + k - \{m(N)\})^2}{(\log N)^3} \right) +\\
      + {}& K\left( \frac{1}{(\log N)^2} - \frac{2(A +k- \{m(N)\})}{(\log N)^3}\right) + 
      \OO\left( \frac{1}{(\log N)^3}\right)
      \end{split}
   \end{equation*}
   So we obtained that
   \begin{equation*}
   \begin{split}
       L=  {}&\log(C_0 p^{-2} q_1 q_2) - E -(k- \{m(N)\}) + \frac{8\log\log N}{C^2(\log N)^2} + 
       \frac{2(C_1-C_0)}{C^2 C_0} \frac{1}{(\log N)^2} + \\
      +{} &\frac{2E}{C \log N} + \frac{2( k- \{m(N)\})}{C \log N} - \frac{4}{C} \frac{(\log\log N)^2}{(\log N)^2} -
       \frac{1}{C^2} \frac{(4 \log \log N)^2}{(\log N)^3} - \frac{C_1-C_0}{C^2 C_0}\frac{4\log \log N}{(\log N)^3} -\\
      - {}& \frac{4}{C} \frac{(\log \log N)( k- \{m(N)\})}{(\log N)^2} - 
      \frac{8}{C^2} \frac{(\log \log N)( k- \{m(N)\})}{(\log N)^3} - \frac{1}{C} \frac{( k- \{m(N)\})^2}{(\log N)^2} + \\
      + {}& \frac{16}{3C} \frac{(\log \log N)^3}{(\log N)^3} + 
      \frac{8}{C} \frac{(\log \log N)^2( k- \{m(N)\})}{(\log N)^3} + 
      \frac{4}{C} \frac{\log \log N( k- \{m(N)\})^2}{(\log N)^3} +\\
      + {}& \frac{C_1-C_0}{C C_0} \left( - \frac{2\log \log N}{(\log N)^2} -
     \frac{4\log \log N}{C(\log N)^3}- \frac{ k- \{m(N)\}}{(\log N)^2} +\frac{4(\log \log N)^2}{(\log N)^3} + \right.\\
      +{}&\left. \frac{4\log \log N (k- \{m(N)\})}{(\log N)^3}\right) + \frac{K}{(\log N)^2} - 
      \frac{4K \log \log N}{(\log N)^3} + \OO\left( \frac{1}{(\log N)^3}\right)
   \end{split}        
   \end{equation*}
   Now let 
\begin{equation*}
    \begin{split}
     E = {}&\frac{8}{C^2}\frac{\log \log N}{(\log N)^2} - \frac{4}{C}\frac{(\log \log N)^2}{ (\log N)^2} - 
     \frac{2(C_1-C_0)}{C C_0} \frac{\log \log N}{(\log N)^2} + \frac{2(C_1-C_0)}{C^2 C_0} \frac{1}{(\log N)^2}- \\
   -{}& \frac{16}{C^2}\frac{(\log \log N)^2}{ (\log N)^3} -\frac{4(C_1-C_0)}{C^2 C_0} \frac{\log \log N}{(\log N)^3} +
    \frac{16}{3C}\frac{(\log \log N)^3}{ (\log N)^3} - \frac{4(C_1-C_0)}{C^2 C_0} \frac{\log \log N}{(\log N)^3} +  \\
     + {}& \frac{K}{(\log N)^2} -\frac{4 K \log \log N}{(\log N)^3} +
     \frac{4(C_1-C_0)}{C C_0} \frac{(\log \log N)^2}{(\log N)^3} + F =\\
     = {}& -\frac{4}{C} \frac{(\log \log N)^2}{(\log N)^2} + \left( \frac{8}{C^2} - 
     \frac{2(C_1-C_0)}{CC_0} \right)\frac{\log \log N}{(\log N)^2} + 
     \left(\frac{2(C_1-C_0)}{C^2C_0} + K\right)\frac{1}{(\log N)^2} + \\
     + {}& \frac{16}{3C} \frac{(\log \log N)^3}{(\log N)^3} + \left( -\frac{16}{C^2} + \frac{4(C_1-C_0)}{C C_0}\right)
     \frac{(\log \log N)^2}{(\log N)^3} - \left(4K + \frac{8(C_1-C_0)}{C^2 C_0}\right)\frac{\log \log N}{(\log N)^3} \\+ {}& F .
    \end{split}
\end{equation*}
To collect those terms which contain $ k- \{m(N)\}$, we introduce the function
\begin{equation*}
    \begin{split}
        H(x) = {}&-x + \frac{2x}{C\log N} - \frac{4}{C}\frac{\log \log N}{(\log N)^2}x - 
        \frac{2}{C}\frac{4\log \log N}{C(\log N)^3}x - \frac{1}{C}\frac{1}{(\log N)^2}x^2 +\\
        + {}& \frac{8(\log \log N)^2}{C(\log N)^3}x + \frac{4(\log \log N)}{C(\log N)^3}x^2 -
        \frac{C_1-C_0}{C C_0} \frac{1}{(\log N)^2}x + \frac{4(C_1-C_0)}{C C_0} \frac{\log \log N}{(\log N)^3}x=\\
        = {}& -x + \frac{2}{C\log N}x - \frac{4}{C}\frac{\log \log N}{(\log N)^2}x - 
        \frac{C_1-C_0}{C C_0} \frac{1}{(\log N)^2}x + \\
        + {}&\left( \frac{4(C_1-C_0)}{C C_0} - \frac{8}{C^2}\right)\frac{\log \log N}{(\log N)^3}x + 
        \frac{8}{C}\frac{(\log \log N)^2}{(\log N)^3}x - \frac{1}{C}\frac{1}{(\log N)^2}x^2 + 
        \frac{4}{C}\frac{\log \log N}{(\log N)^3}x^2 .
    \end{split}
\end{equation*}
Inserting these expressions, we obtain
\begin{equation*}
    \begin{split}
        L = {}&\log(C_0 p^{-2}q_1q_2) + H( k- \{m(N)\}) - F + \frac{2E}{C \log N} + O\left( \frac{1}{(\log N)^3}\right)
        =\\
        = {}&\log(C_0 p^{-2}q_1q_2) + H( k- \{m(N)\}) - F + \frac{16}{C^3}\frac{\log \log N}{(\log N)^3} - 
        \frac{8}{C^2}\frac{(\log \log N)^2}{(\log N)^3} - \\
        - {}& \frac{4(C_1-C_0)}{C^2 C_0}\frac{\log \log N}{(\log N)^3} + \OO\left( \frac{1}{(\log N)^3}\right) .
    \end{split}
\end{equation*}
So choosing
\begin{equation*}
    \begin{split}
        F = \frac{16}{C^3}\frac{\log\log N}{(\log N)^3} - \frac{8}{C^2}\frac{(\log\log N)^2}{(\log N)^3} - 
        \frac{4(C_1-C_0)}{C^2 C_0} \frac{\log \log N}{(\log N)^3} ,
    \end{split}
\end{equation*}
we have
\begin{equation*}
    L = \log(C_0 p^{-2}q_1q_2) + H( k- \{m(N)\}) + \OO\left( \frac{1}{(\log N)^3}\right).
\end{equation*}
   So we obtain that
   \begin{equation*}
       \begin{split}
           P(\mu(N)- [m(N)] < k) & = e^{-l} \left(1+ \OO\left(\frac{(\log N)^3}{N}\right) \right) \\        
           &=e^{-(1/p)^L} \left(1+ \OO\left(\frac{(\log N)^3}{N}\right) \right) \\
           &=e^{-(1/p)^{\log(C_0 p^{-2}q_1q_2) + H( k- \{m(N)\}) + \OO\left( \frac{1}{(\log N)^3}\right)} }
           \left(1+ \OO\left(\frac{(\log N)^3}{N}\right) \right) \\         
           &=\exp\left(-p^{-(\log(C_0 p^{-2}q_1q_2) + H( k- \{m(N)\}) )}\right)
           \left(1+ \OO\left(\frac{1}{(\log N)^3}\right) \right) .          
       \end{split}
   \end{equation*}
    \end{proof}
\section{Simulation results} \label{simu}
Analysing the beginning of the proof of Theorem \ref{T2}, we can see that the lemma of \citet{csaki1987limit}
offers good approximation if $p$ is small, but it does not offer good approximation if $p$ is close to $1$.
However, our simulation studies show, that the approximation for the longest run is very good for small values of $p$,
but it is still appropriate if $p$ is close to $1$.

We performed several computer simulation studies for certain fixed values of
$p$, $q_1$, and $q_2$.
Below $N$ denotes the length of the sequences generated by us and $s$ denotes the number of repetitions 
on the $N$-length sequences.

\begin{figure}[h!]
    \centering
    \begin{subfigure}[b]{0.45\textwidth}
        \includegraphics[width=\textwidth]{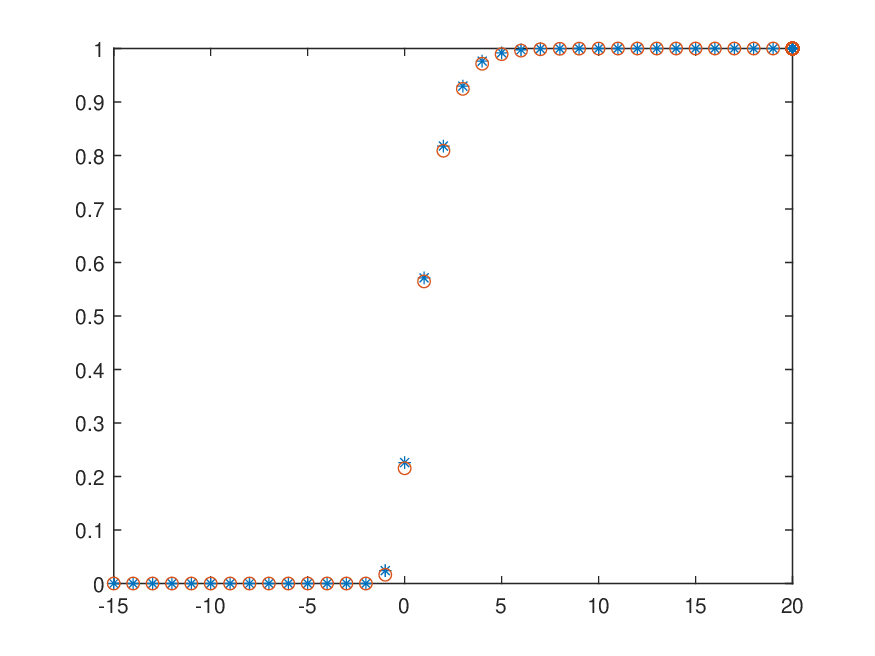}
        \caption*{(a) Longest run}
    \end{subfigure}
    \begin{subfigure}[b]{0.45\textwidth}
        \includegraphics[width=\textwidth]{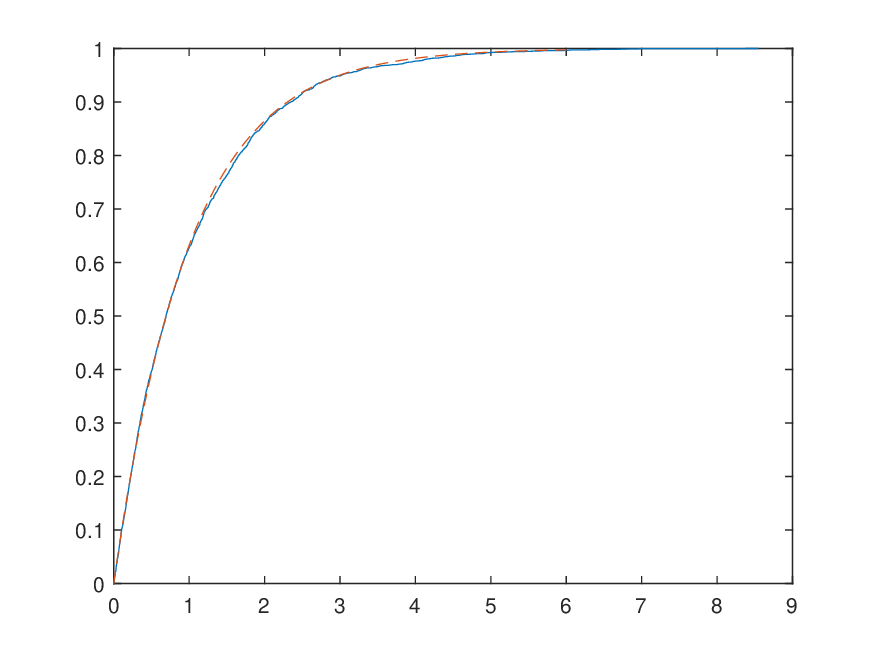}
        \caption*{(b) First hitting time, $m=16$}
    \end{subfigure}   
    \caption{Longest at most $1+1$ contaminated run and the first hitting time 
    when  $p=1/3$, $q_1=1/3$, $q_2=1/3$, $N= 3 \cdot 10^6$, $s=3000$}  \label{TT333run}
 \end{figure}

\begin{figure}[h!]
    \centering
    \begin{subfigure}[b]{0.45\textwidth}
        \includegraphics[width=\textwidth]{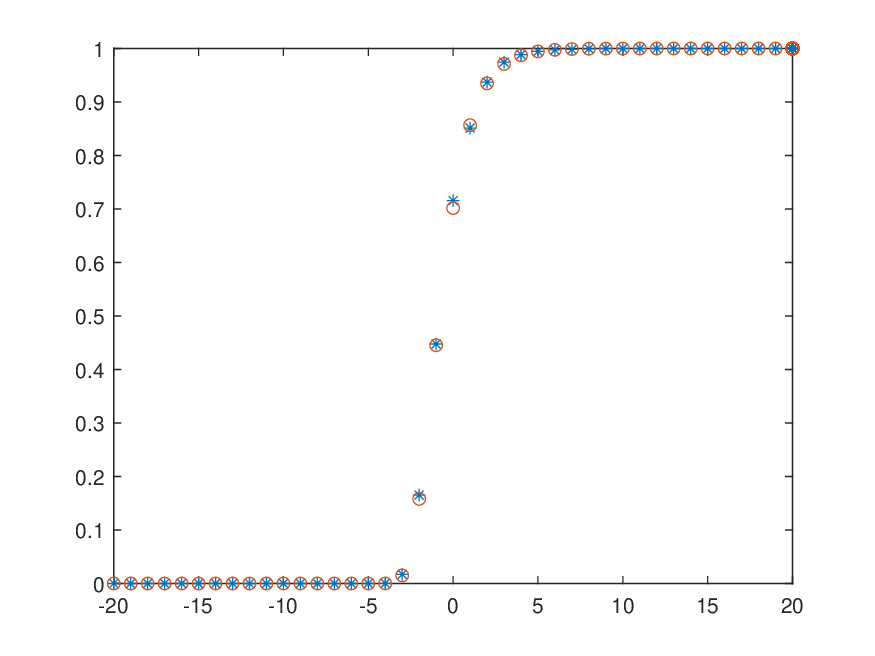}
        \caption*{(a) Longest run}
    \end{subfigure}
    \begin{subfigure}[b]{0.45\textwidth}
        \includegraphics[width=\textwidth]{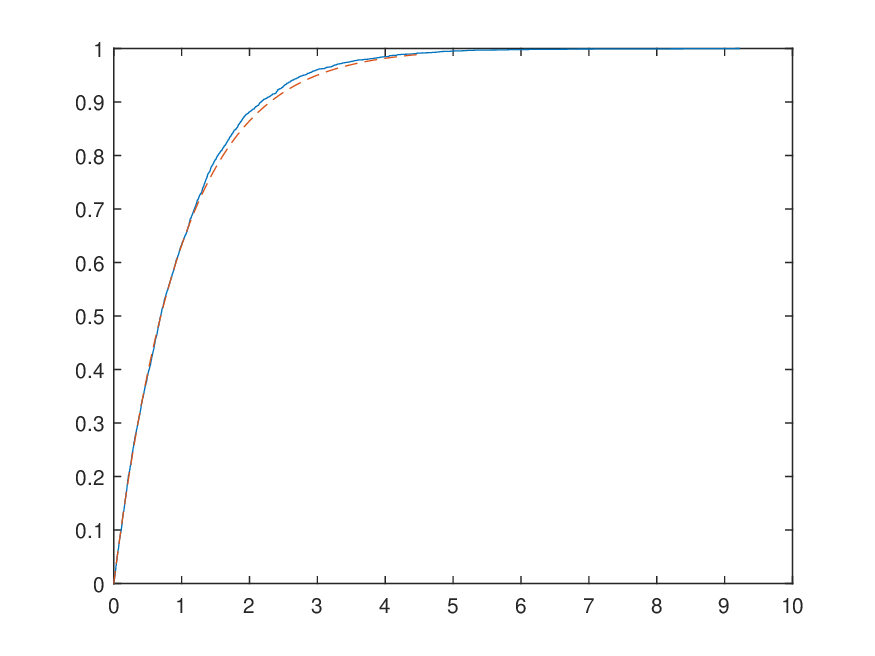}
        \caption*{(b) First hitting time, $m=19$}
    \end{subfigure}   
    \caption{Longest at most $1+1$ contaminated run and the first hitting time 
    when  $p=0.4$, $q_1=0.3$, $q_2=0.3$, $N= 3 \cdot 10^6$, $s=3000$}  \label{TT04run}
 \end{figure}

\begin{figure}[h!]
    \centering
    \begin{subfigure}[b]{0.45\textwidth}
        \includegraphics[width=\textwidth]{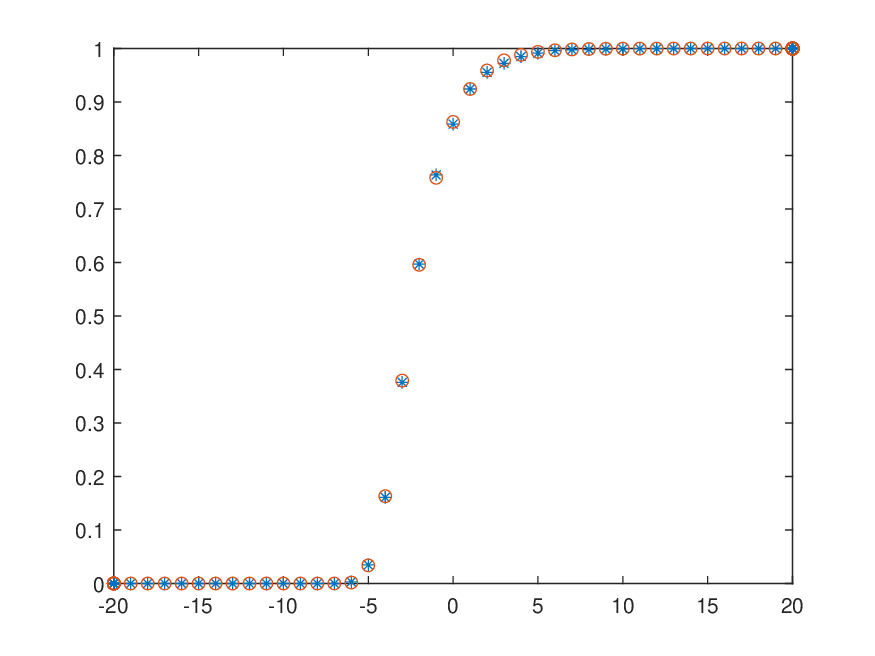}
        \caption*{(a) Longest run}
    \end{subfigure}
    \begin{subfigure}[b]{0.45\textwidth}
        \includegraphics[width=\textwidth]{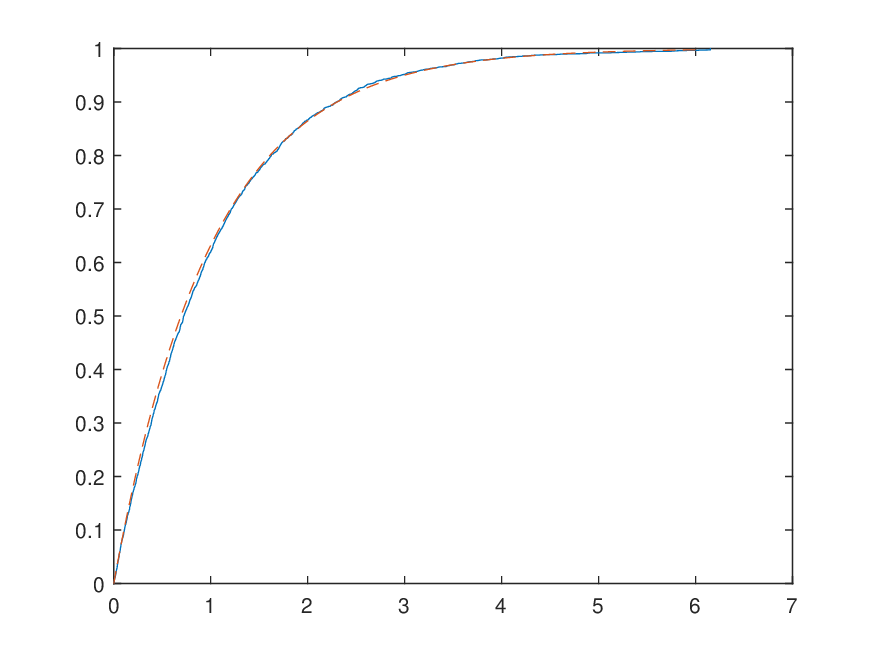}
        \caption*{(b) First hitting time, $m=25$}
    \end{subfigure}   
    \caption{Longest at most $1+1$ contaminated run and the first hitting time 
    when  $p=0.5$, $q_1=0.4$, $q_2=0.1$, $N= 4 \cdot 10^6$, $s=3000$}  \label{TT541run}
 \end{figure}

\begin{figure}[h!]
    \centering
    \begin{subfigure}[b]{0.45\textwidth}
        \includegraphics[width=\textwidth]{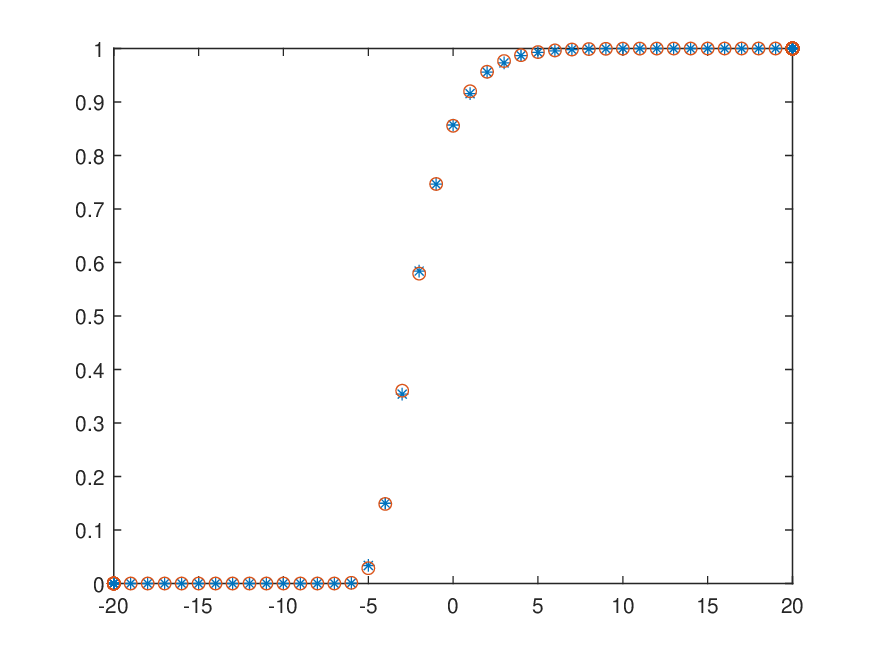}
        \caption*{(a) Longest run}
    \end{subfigure}
    \begin{subfigure}[b]{0.45\textwidth}
        \includegraphics[width=\textwidth]{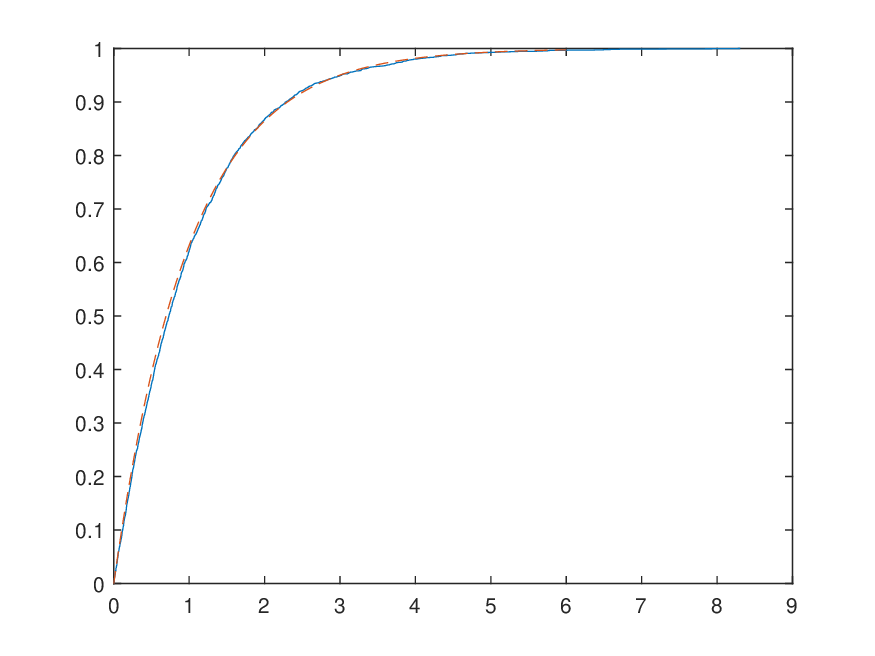}
        \caption*{(b) First hitting time, $m=23$}
    \end{subfigure}   
    \caption{Longest at most $1+1$ contaminated run and the first hitting time 
    when  $p=0.5$, $q_1=0.3$, $q_2=0.2$, $N= 3 \cdot 10^6$, $s=3000$}  \label{TT532run}
 \end{figure}

\begin{figure}[h!]
    \centering
    \begin{subfigure}[b]{0.45\textwidth}
        \includegraphics[width=\textwidth]{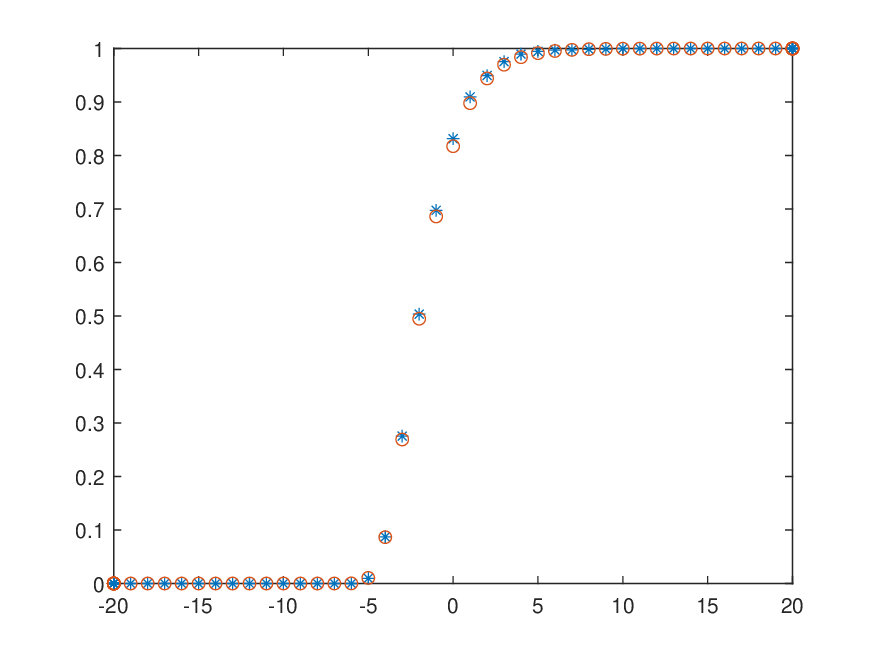}
        \caption*{(a) Longest run}
    \end{subfigure}
    \begin{subfigure}[b]{0.45\textwidth}
        \includegraphics[width=\textwidth]{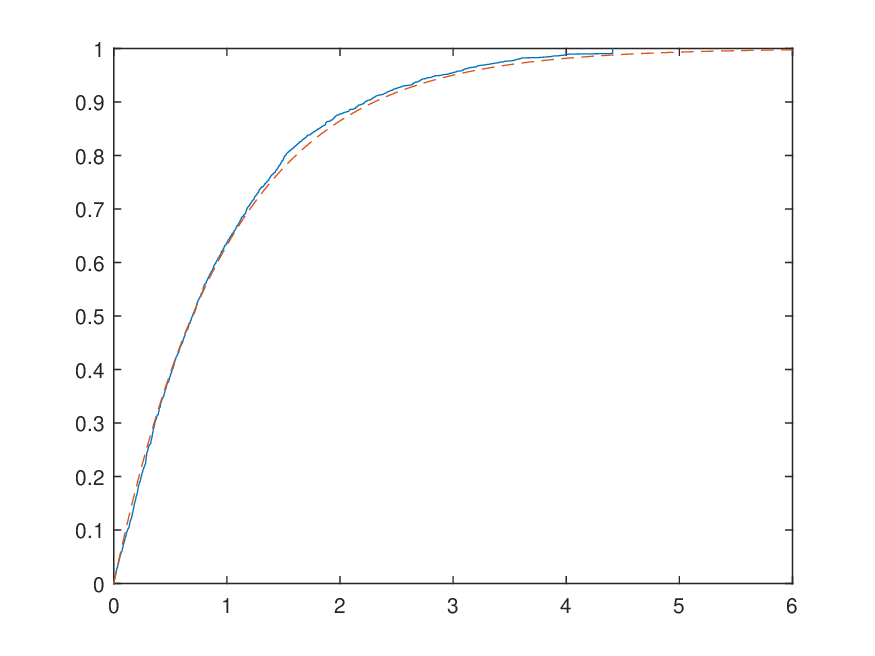}
        \caption*{(b) First hitting time, $m=25$}
    \end{subfigure}   
    \caption{Longest at most $1+1$ contaminated run and the first hitting time 
    when  $p=0.5$, $q_1=0.25$, $q_2=0.25$, $N= 2 \cdot 10^6$, $s=2000$}  \label{TT05run}
 \end{figure}

\begin{figure}[h!]
    \centering
    \begin{subfigure}[b]{0.45\textwidth}
        \includegraphics[width=\textwidth]{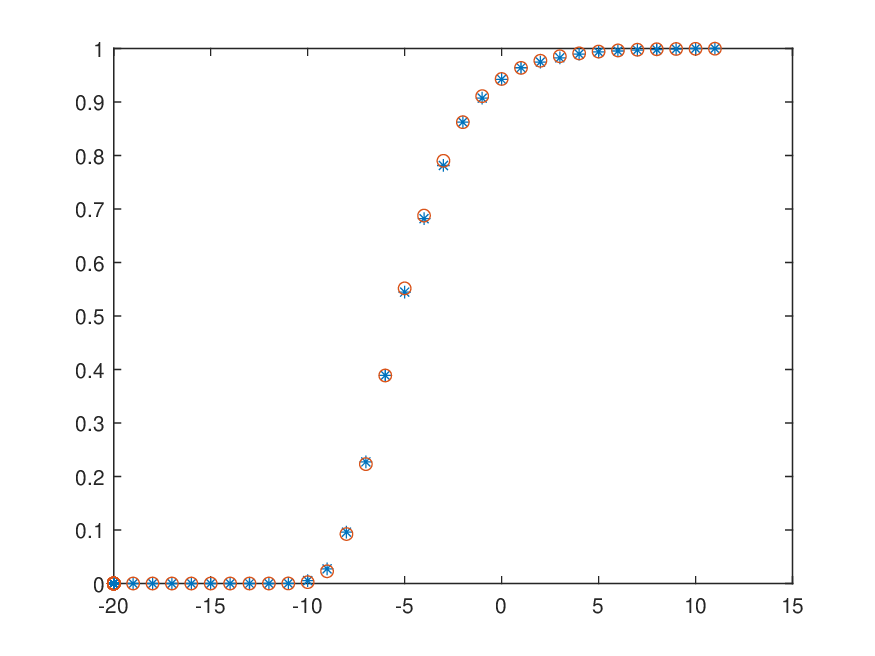}
        \caption*{(a) Longest run}
    \end{subfigure}
    \begin{subfigure}[b]{0.45\textwidth}
        \includegraphics[width=\textwidth]{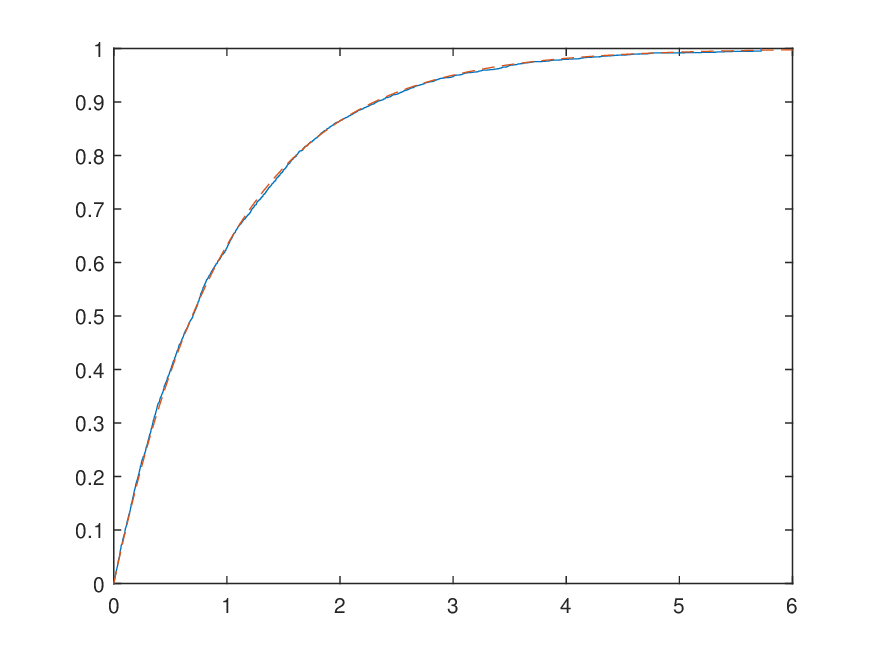}
        \caption*{(b) First hitting time, $m=34$}
    \end{subfigure}   
    \caption{Longest at most $1+1$ contaminated run and the first hitting time 
    when  $p=0.6$, $q_1=0.2$, $q_2=0.2$, $N= 4 \cdot 10^6$, $s=3000$}  \label{TT06run}
 \end{figure}

\begin{figure}[h!]
    \centering
    \begin{subfigure}[b]{0.45\textwidth}
        \includegraphics[width=\textwidth]{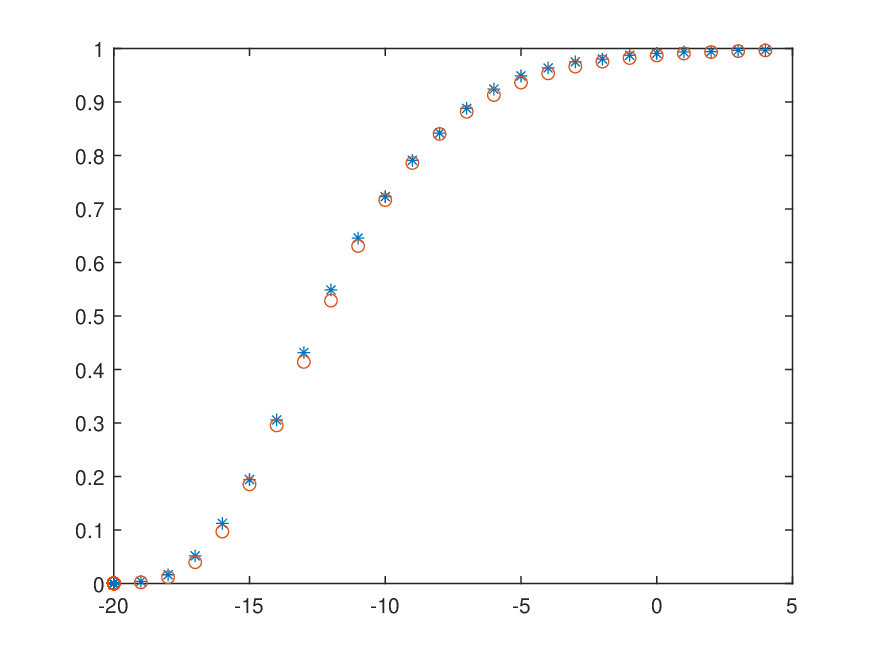}
        \caption*{(a) Longest run}
    \end{subfigure}
    \begin{subfigure}[b]{0.45\textwidth}
        \includegraphics[width=\textwidth]{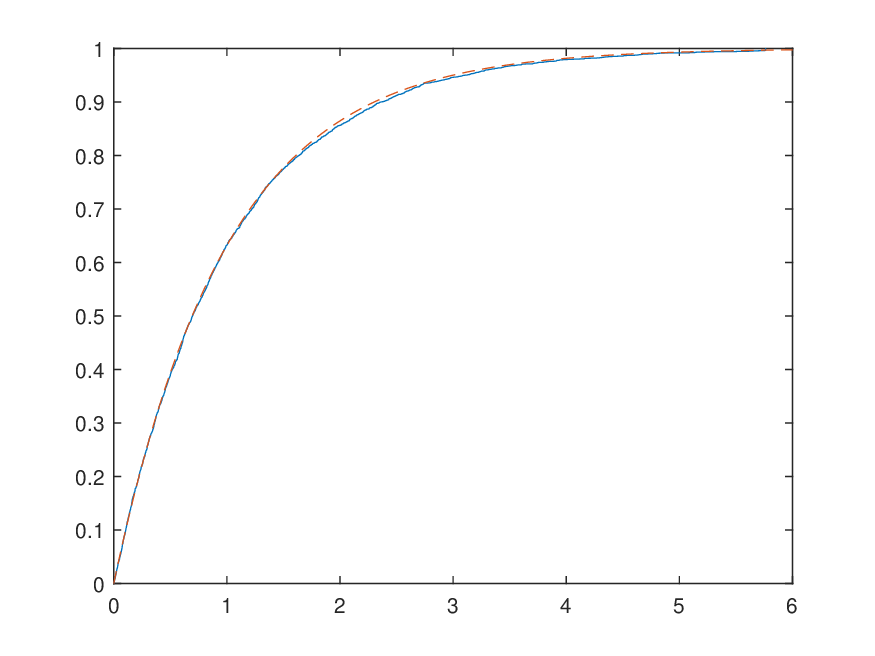}
        \caption*{(b) First hitting time, $m=47$}
    \end{subfigure}   
    \caption{Longest at most $1+1$ contaminated run and the first hitting time 
    when  $p=0.7$, $q_1=0.2$, $q_2=0.1$, $N= 4 \cdot 10^6$, $s=3000$}  \label{TT07run}
 \end{figure}

\begin{figure}[h!]
    \centering
    \begin{subfigure}[b]{0.45\textwidth}
        \includegraphics[width=\textwidth]{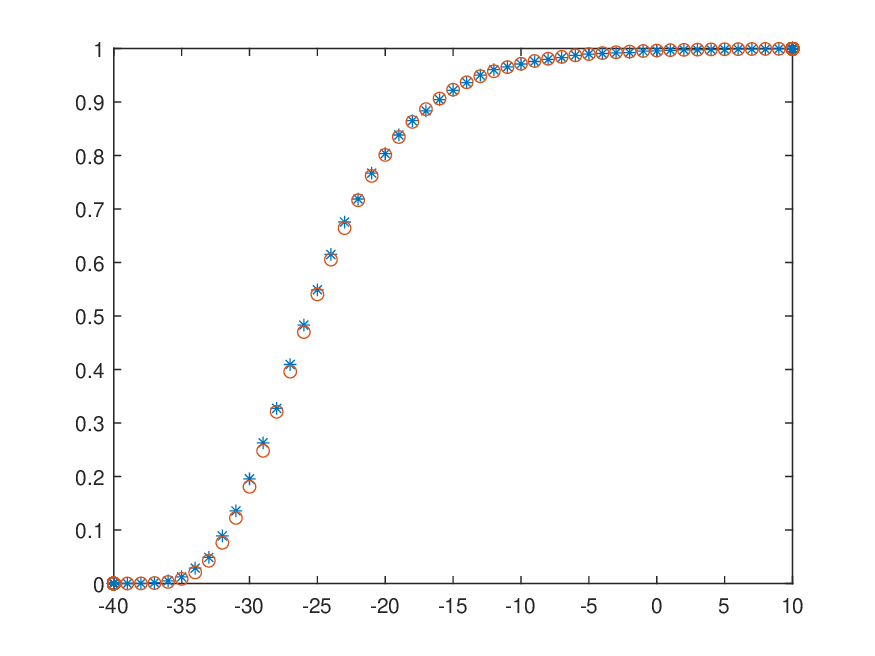}
        \caption*{(a) Longest run}
    \end{subfigure}
    \begin{subfigure}[b]{0.45\textwidth}
        \includegraphics[width=\textwidth]{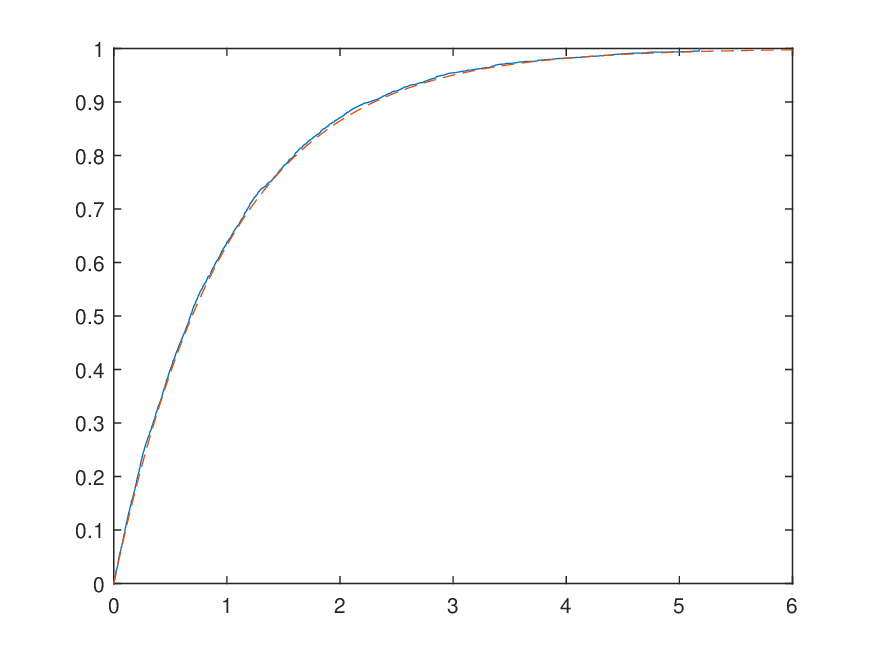}
        \caption*{(b) First hitting time, $m=72$}
    \end{subfigure}   
    \caption{Longest at most $1+1$ contaminated run and the first hitting time 
    when  $p=0.8$, $q_1=0.1$, $q_2=0.1$, $N= 3 \cdot 10^6$, $s=3000$}  \label{TT08run}
 \end{figure}

Figures \ref{TT333run}- \ref{TT08run} present the results of the simulations.
The left hand side part of each figure shows the empirical distribution of the longest at most $1+1$ contaminated run
and its approximation suggested by Theorem \ref{T2}.
Asterisk (i.e. $\ast$) denotes the result of the simulation, i.e. the  empirical distribution of 
the longest at most $1+1$ contaminated run and circle ($\circ$) denotes approximation offered by Theorem \ref{T2}.
The right hand side of each figure shows the first hitting time of the $m$-length at most $1+1$ contaminated run.
Solid line shows the result of the simulation for the distribution function and dashed line shows
the distribution function $1-e^{-x}$ suggested by our Theorem \ref{T1}.

\section{Appendix} \label{app}    
\setcounter{equation}{0}
The following lemma of Cs\'aki, F\"oldes and Koml\'os plays a fundamental role in our proofs.
\begin{lemma}	\label{MLsf}
	(Main lemma, stationary case, finite form of \citet{csaki1987limit}.)
	Let $X_1,X_2,\dots$ be any sequence of independent random variables, 
and let $\mathcal{F}_{n,m}$ be the $\sigma$-algebra generated by the random variables $X_n,X_{n+1},\dots,X_{n+m-1}$. 
Let $m$ be fixed and let $A_n=A_{n,m}\in\mathcal{F}_{n,m}$. 
Assume that the sequence of events $A_n = A_{n,m}$ is stationary, that is
$P(A_{i_1+d} A_{i_2+d} \cdots A_{i_k+d})$ is independent of $d$.
	 
	Assume that there is a fixed  number $\alpha$, $0<\alpha \le 1$, 
	such that the following three conditions hold for some fixed $k$ with $2 \le k \le m$, 
	and fixed $\varepsilon$ with $0< \varepsilon < \min\{ p/10, 1/42\}$ 
	\begin{enumerate}[(SI)]
		\item 
		\begin{equation*}
			|P(\bar{A_2} \cdots \bar{A_k}| A_1)-\alpha| <\varepsilon ,
		\end{equation*} 
		\item
		\begin{equation*}
			\sum_{k+1 \le i \le 2m} P(A_i|A_1) <\varepsilon ,
		\end{equation*}
		\item
		\begin{equation*}
			P(A_1) < \varepsilon/m .
		\end{equation*}
	\end{enumerate}
	Then, for all $N > 1$,
	\begin{equation*}
		\left|\frac{P(\bar{A_2} \cdots \bar{A}_N|A_1)}{P(\bar{A_2} \cdots \bar{A}_N)} - \alpha \right| <7 \varepsilon
	\end{equation*}
	and 
	\begin{equation} \label{MLsfe}
		e^{-(\alpha+10\varepsilon)NP(A_1)-2mP(A_1)} < P(\bar{A_1} \cdots \bar{A}_N) <e^{-(\alpha-10\varepsilon)NP(A_1)+2mP(A_1)}.
	\end{equation}
\end{lemma}
   \bibliographystyle{apalike}
\bibliography{TwoTypNew+}
\end{document}